\documentclass[A4]{article}
\usepackage{amssymb,amsmath,latexsym, amscd, epsfig}
\newtheorem{theorem}{Theorem}
 
\newtheorem{corollary}[theorem]{Corollary}

\newtheorem{remark}{Remark:} 

\newcommand{\Z}{\mathbb Z}
\newcommand{\F}{\mathcal F}

\begin{document}
\title{On a theorem of Goussarov}
\author{Jim Conant \footnote{Supported by the NSF through a VIGRE postdoctoral position, and by Max-Planck-Institut where the author was a guest when he
wrote this paper}}
\date{\today}
\maketitle
\begin{abstract}
In this paper, the easier methods of my thesis are applied to give a simple
proof of a theorem of Goussarov. This theorem relates two possible notions
of finite type equivalence of knots, links, or string links,
 showing that the resulting filtrations are the same up to a degree shift by
a factor of two. This is then applied to the situation of rooted claspers to 
show that rooted clasper surgeries of sufficiently high degree must 
preserve type $k$ invariants. As a consequence, grope cobordisms of
sufficiently high class must preserve type $k$ invariants.
This result is applied in \cite{ct} to show Theorem 2 of that 
paper.
\end{abstract}
\section{Introduction}
The heart of this paper is a simple proof of a theorem of Goussarov,
using an elementary tool developed in my thesis. His theorem relates two 
possible notions of finite type invariants, the standard one which he
co-invented with Vassiliev, and a more subtle one based on what he calls
``interdependent'' modifications of links. He proves that these invariants
coincide up to a degree shift by a factor of two. 

The reason I wrote this paper is Corollary \ref{cor}, that a grope cobordism
preserves Vassiliev invariants up to roughly half the degree.
 This corollary is needed to prove
Theorem 2 of \cite{ct}, which states that grope cobordism of class $c$
is precisely the same as
surgery on \emph{simple} claspers (see \cite{h} or \cite{ct}) with 
\emph{grope degree} $c$ or more. It is much easier to see that 
 grope cobordism coincides
with \emph{rooted} clasper surgeries, for which Goussarov's 
interdependent modifications are well-suited.

In my thesis I considered a problem similar to Corollary \ref{cor},
 the difference being that I considered
knots bounding gropes, which is a much stronger notion than cobounding a
grope with the unknot. In that case, I obtained a sharp answer, namely
a knot bounding a grope of class $k$ is $\lceil k/2\rceil$-trivial. The 
current paper implies $\lfloor (k-1)/2 \rfloor$, which is not very far
away. Achieving the sharp result is a lot of work for not much gain. 

Goussarov's main idea was to define filtrations using alternating sums over
``disjoint variations.''
For instance, one may take the variations to be crossing changes of a knot,
as he did in \cite{g1}, defining the usual finite type filtration 
independently of Vassiliev. In \cite{g2} he considered a more subtle notion
of variation, an ``interdependent modification,''
which involves replacing an arc of a circle with an arbitrary
other arc in the link complement with the same endpoints. Evidently he was 
trying to find more subtle invariants than those which are of standard finite
type. But his discovery, which is the main subject of this paper, was
that they are the same up to a degree change.
 Most recently, he invented
and developed
a notion of 3-manifold finite type invariants, where the variation is a surgery
along a ``Y.'' (Equivalently, a Matveev move on a genus 3 handlebody.)
This notion appears in \cite{ggp} for example.
In \cite{gr} Garoufalidis and Rozansky use a variation of this notion to study
pairs (Homology 3-sphere, Knot), where the leaves of the ``Y'' must link
the knot trivially. 
If it weren't for his untimely death, Goussarov surely would
have produced more exotic and beautiful notions.

Acknowledgments: It is a pleasure to thank Thang Le and Peter Teichner for
useful discussions. I also wish to dedicate this to the late Mikhail
Goussarov, whose paper \cite{g1}, gave me an inspiring start in this field. 

\section{Two filtrations}
Let $X$ be the set of isotopy classes of links or string links with a 
fixed number of components. The Vassiliev filtration
of $\Z[X]$ is a descending filtration $\Z[X]=\F^v_0\supseteq \F^v_1 \supseteq 
\F^v_2\supseteq \cdots$. Each piece $\F^v_k$ is spanned by alternating sums
of the following form. Let $x\in X$ be a knot, link, or string link. Choose
$k+1$ framed arcs from $x$ to itself, which guide homotopies of $x$ supported in
neighborhoods of these arcs, which take the little piece of $x$
at one end of the arc and push this across the other end. Such a move is 
called a finger move. (Such finger moves can always be realized as crossing
changes in some projection of $x$.) Let the set of these finger moves be called
$S$. Then define $[x;S]= \sum_{\sigma\subset S}(-1)^{|\sigma|}x_\sigma$ where
$x_\sigma$ is $x$ modified by the finger moves in $\sigma$. Now, by definition
$\F^v_k$ is additively generated by elements of the form $[x,S]$, where $S$ 
has $k+1$ finger moves in it.

One can replace single finger moves by groups\footnote{I use the term group
to indicate that a collection of moves is to be viewed as a single variation.}
 of finger moves in the above 
definition, and this is well-known to give the same filtration. See \cite{g1}
for a proof.

An alternative filtration can be defined in the same way using moves more 
general than finger moves. One way to think of these moves is that they are
guided by circles which are attached to $x$ along a subarc. 
The move is to replace the part of $x$ running across the circle with 
the other arc of the circle. $k+1$ disjoint circles lead to generators of this 
alternative filtration, which I denote by $\F^{alt}_k$. Goussarov calls these
interdependent moves. The reason is that such moves cannot neccesarily
be realized by independent, i.e. disjoint, groups of finger moves.

As before, one can replace single moves along circles by groups of moves along
circles and the same filtration is achieved, the proof being modified
\emph{mutatis mutandis}.

The main result of \cite{g2}, which is stated there in the dual setting,
 is that $\F^v_{k} = \F^{alt}_{2k}=\F^{alt}_{2k+1}$.
It is proved in two steps. 

\begin{theorem}\label{step1}
$\F^v_k\subset \F^{alt}_{2k+1}$.
\end{theorem}

\begin{theorem}\label{step2}
 $\F^{alt}_{2k}\subset \F^v_k$.
\end{theorem}

Theorem \ref{step1} is the easier of the two. I provide a repackaged 
version of Goussarov's proof for completeness. Theorem \ref{step2} is harder 
and this is the one for which I provide a new simplified proof. 

\section{Proofs of Theorems \ref{step1} and \ref{step2}}

\noindent \emph{[Proof of Theorem \ref{step1}]}

 Each finger move has two interdependent moves 
associated to it as follows. A symmetric way to view a finger move 
along an arc $\eta$
is
to push the strands of $x$ at both endpoints of $\eta$ along $\eta$ so
that they crash through each other at some point in the middle of $\eta$.
For each of the two modified strands of $x$, there is a disk cobounding
that strand with the strand before the finger move, the disk having
been swept out by the isotopy. These two disks intersect in a single clasp.
The boundaries of these two disks are the circles guiding the two
interdependent moves associated to the finger move.  
Now I claim the alternating sum over a set of finger moves is, up
to a sign, the same
as the alternating sum over the associated interdependent moves. This is 
essentially because
if one only does one of a given pair associated to a finger move, this is 
just an isotopy. Let $S=\{s_i\}$ be the set of $k+1$ finger moves,
and $T=\{a_i,b_i\}$ the associated set of $2k+2$ interdependent moves. 
Inductively I show that 
$$\sum_{\sigma\subset T} (-1)^{|\sigma |}x_\sigma = (-1)^{k}
\sum_{\mu\subset S} (-1)^{|\mu |}x_\mu,$$
which will be sufficient to prove the theorem.
This left hand sum breaks up into
$$\sum_{\tau\subset\{a_2,\ldots, b_{k+1}\}}(-1)^{|\tau|}(x_\tau - 
x_{\tau\cup\{a_1\}} - x_{\tau\cup\{b_1\}} + x_{\tau\cup\{a_1,b_1\}}),$$ and since
the first three terms are equal, we get
$$-\sum_{\tau\subset\{a_2,\ldots, b_{k+1}\}}(-1)^{|\tau|}x_\tau + 
\sum_{\tau\subset\{a_2,\ldots, b_{k+1}\}}(-1)^{|\tau |}(x_{s_1})_\tau,$$ and
by induction this is just
$$(-1)^{k}\sum_{\mu\subset \{s_2,\ldots,s_{k+1}\}}((-1)^{|\mu |}x_\mu 
- (-1)^{|\mu |}x_{\{s_1\}\cup \mu})$$
which equals $(-1)^k\sum_{\mu\subset S}(-1)^{|\mu|}x_\mu$ as desired.
$\hfill\Box$

\vspace{1em}
\noindent \emph{[Proof of Theorem \ref{step2}]}

Given a set of $2k+1$ interdependent moves, $S$ on $x$, we wish to show
that $[x;S]\in \F^v_{k}$. 
The strategy is to show that $[x;S]$ is congruent modulo $\F^v_k$ to
sums of simpler alternating sums $[x;S^\prime]$, where the complexity
of such alternating sums
is recorded by a graph I will define in a moment. After iteration, we reduce
the problem to alternating sums for which a direct argument is possible
to show that they are congruent to $0$.

Fix a projection of $x$ and the circles in $S$. 
(That is, think of it as a planar picture with over and under crossing
data.) 
The projection is assumed to be such that $x$ does not cross itself
where circles are attached.
In order to keep track of the complexity of the interdependence
of the moves, I now define a graph associated to the picture. It has
$2k+1$ vertices corresponding to the $2k+1$ circles that guide the 
interdependent moves. Fix an ordering on these vertices, say they are
called $v_1,\ldots, v_{2k+1}$. Draw an edge between $v_i$ and $v_j$ if $i<j$,
but the circle corresponding to $v_i$ crosses over the one corresponding
to $v_j$. Draw an edge from a vertex to itself if the corresponding
circle is knotted. Finally, label a vertex with a star if $x$ crosses over
the corresponding circle. 

For each edge of the graph, one can do an obvious group of crossing changes
of $x$ union the circles, which eliminates that edge. That is, to eliminate an edge between $v_i$ and $v_j$, where $i<j$, do those crossing changes which
make $v_i$ always pass under $v_j$ at each crossing. Similarly, to eliminate
a star, there is a group of crossing changes that always make $x$ pass
under a given circle.
Finally there is a
group of crossing changes eliminating a self-edge, i.e. by unknotting the
corresponding circle.

Suppose a graph has at least $k+1$ edges plus stars.
Let $T$ denote a set of $k+1$ groups of crossing changes on $x$ union the
circles, each of which removes an edge or star from the graph.
Each such group of crossing changes  induces a group of
crossing changes on each summand $x_\sigma$ of $[x;S]$.
Thus each $x_\sigma$ is congruent modulo $\F^v_k$ to 
$\sum_{\emptyset\neq\tau\subset T} \pm {x_\sigma}_\tau$. 
Then 
\begin{gather*}
\begin{split}
[x;S] &= \sum_{\sigma\subset S}(-1)^{|\sigma |}x_\sigma\\ 
&\equiv \sum_{\sigma\subset S}(-1)^{|\sigma |}\sum_{\emptyset\neq\tau\subset T}
\pm (x_\sigma)_\tau\\
&= \sum_{\emptyset\neq\tau\subset T}\pm\sum_{\overline{\sigma}\subset S_\tau}
(-1)^{|\overline{\sigma} |} x_{\overline{\sigma}}\\ 
&=\sum_{\emptyset\neq\tau\subset T}\pm [x,S_\tau],
\end{split}
\end{gather*}
where $S_\tau$
is the set of interdependent moves guided by the circles modified by $\tau$.
Thus it suffices to show that $[x;S_\tau]\in\F^v_k$ for 
$\tau\neq\emptyset$.
Each of the resulting sets $S_\tau, \tau\neq\emptyset$,
of interdependent moves has fewer edges plus stars
in the resulting graph.
  We can always iterate this simplification unless
the number of edges plus stars is less than or equal to $k$. It is hence 
sufficient to consider this case.
 
On such a graph, I claim that there are $k+1$ unstarred vertices, no pair
of which has a connecting edge. This follows from the claim that there
are at least $k+1$ connected components without any stars on them. 
To see this, let $st$ denote the number of stars. Then $E+st\leq k$ implies
$E\leq k-st$. The Euler characteristic can be computed in two
different ways: $V-E = b_0-b_1$, so that
$b_0\geq b_0-b_1 = 2K+1-E \geq (k+1)+st$.
Thus
$b_0 - st \geq k+1$.

The fact that the vertices are unstarred means they each bound a disk, and
the lack of edges implies the disks are disjoint. There are now $k+1$ groups of
finger
moves of $x$ union the circles which push everything out of the $k+1$ disks.
Thus, as above, modulo $\F^v_k$, we can assume that at least one of the circles
guiding a move bounds an embedded disk with no intersections with anything
in its interior. This means that the move is an isotopy. However
$[x;S]$ is obviously $0$ if any of the moves in $S$ are isotopies. $\hfill\Box$

\section{Application to rooted claspers and grope cobordism}
 For the reader's 
convenience in following the proof of Theorem \ref{grope}, I remind him that
a rooted clasper is a clasper in the sense of Habiro(see \cite{h}), embedded
in the complement of a knot or link, such that there is one root leaf, which is
a zero framed leaf linking the knot as a little meridian, whereas the 
other leaves can be embedded arbitrarily.

\begin{theorem}\label{grope}
Suppose two knots $K_1$ and $K_2$ are related by a rooted clasper surgery
 with $2k+1$ non-root leaves. Then $K_1-K_2\in \F^v_k = \F^{alt}_{2k}$.
\end{theorem}

\noindent\emph{[Proof]}

 Surgering along the clasper, the original knot
is modified inside a regular neighborhood of the clasper union the disk
bounding its root leaf to get $K_2$. We will find $2k+1$ groups of interdependent moves as follows.
For each non-root leaf of the clasper there is an interdependent move of a 
subarc of the leaf so that the result is a tiny leaf bounding a disk avoiding
intersections with anything. Each of these moves on the clasper descend to 
a group of interdependent moves
which take strands of the knot running through the neighborhood of the leaf
and move them to a position corresponding to the modified clasper. Thus we 
have found a set $S$ of $2k+1$ groups of interdependent moves on $K_2$.
Thus the surgered knot $K_2$, is congruent
modulo $\F^{alt}_{2k}$ to the sum $-\sum_{\emptyset\neq\sigma\subset S}
(-1)^{|\sigma|}(K_2)_\sigma$, but for each nonempty $\sigma$, $(K_2)_\sigma$
is by construction just $K_1$ modified by a clasper with at least one
trivial leaf. Such clasper surgeries do not change the knot, hence
$K_2$ is congruent to $-\sum_{\emptyset\neq\sigma\subset S}  (-1)^{|\sigma|} 
K_1 = K_1$.$\hfill\Box$
 
\begin{corollary}\label{cor}
Suppose two knots $K_1$ and $K_2$ are related by a grope
cobordism
 of class $2k+1$. Then $K_1-K_2\in \F^v_k$.
\end{corollary}

\noindent\emph{[Proof]}

By Proposition 6 of \cite{ct}, we may assume that the grope cobordism has 
all of its stages of genus one. By Theorem 9 of the same paper, this is just
a rooted clasper surgery with $2k+1$ non-root leaves. 
$\hfill\Box$

\begin{remark}
Even though this Corollary is used in \cite{ct} to prove Theorem 2 of that 
paper, Theorem 2 is not used to prove Theorem 9 or Proposition 6, so there is 
no logical circuity. Corollary \ref{cor} could have been
proven directly in the same way as Theorem \ref{grope}, using techniques
of my thesis \cite{c}, but as there is no logical neccessity I have
avoided the added complication.
\end{remark}

\end{document}